\documentclass{amsart}
\usepackage{amsfonts}
\usepackage{amssymb}
\usepackage{amsmath}
\usepackage{latexsym}
\usepackage{hyperref}

\newcommand{\be}{\begin{equation}}
\newcommand{\ee}{\end{equation}}
\newcommand{\bea}{\begin{eqnarray}}
\newcommand{\eea}{\end{eqnarray}}

\begin{document}
\title{Weak Ricci curvature bounds for Ricci shrinkers}
\author{Bennett Chow}
\author{Peng Lu}
\author{Bo Yang$^{1}$}

\begin{abstract}
We show that for a complete Ricci shrinker there exists a sequence of
points tending to infinity whose norms of the Ricci tensor grow at
most linearly.
\end{abstract}
\maketitle

All objects are $C^{\infty }$. Let $\left( \mathcal{M}^{n}\!,g\right) $ be a
Riemannian manifold and $\phi ,f:\mathcal{M}\rightarrow \mathbb{R}$.\ For $%
\gamma :\left[ 0,\bar{s}\right] \rightarrow \mathcal{M}$, $\bar{s}>0$,
define $S=\gamma ^{\prime }$ and $\mathcal{J}\left( \gamma \right)
=\int_{0}^{\bar{s}}\left( \left\vert S\left( s\right) \right\vert ^{2}+2\phi
\left( \gamma \left( s\right) \right) \right) ds$. A critical point $\gamma $
of $\mathcal{J}$ on paths with fixed endpoints, called a $\phi $-geodesic,
satisfies $\nabla _{S}S=\nabla \phi $ and $\left\vert S\right\vert
^{2}-2\phi =C$. Let $\operatorname{Rc}_{f}=\operatorname{Rc}+\nabla \nabla f$. For a minimal
$\phi $-geodesic,\vspace{-0.04in}%
\begin{equation}
-\int_{0}^{\bar{s}}\zeta ^{2}\Delta _{f}\phi ds+\int_{0}^{\bar{s}}\zeta ^{2}%
\operatorname{Rc}_{f}\left( S,S\right) ds\leq \int_{0}^{\bar{s}}\left( n\left( \zeta
^{\prime }\right) ^{2}-2\zeta \zeta ^{\prime }\left\langle \nabla
f,S\right\rangle \right) ds,\vspace{-0.04in}  \label{second variation}
\end{equation}%
where $\Delta _{f}=\Delta -\nabla f\cdot \nabla $ and $\zeta :\left[ 0,\bar{s%
}\right] \rightarrow \mathbb{R}$ is piecewise $C^{\infty }$, vanishing at $0$
and $\bar{s}$.

Let $\left( g,f\right) $ be a complete shrinker and satisfy $\operatorname{Rc}_{f}=%
\frac{1}{2}g$ and $f-\left\vert \nabla f\right\vert ^{2}=R>0$. Let $c>0$ and
$2\phi =c\frac{R}{f}$. From $\Delta _{f}R=-2\left\vert \operatorname{Rc}\right\vert
^{2}+R$ and $\Delta _{f}f=\frac{n}{2}-f$ we compute\vspace{-0.04in}%
\begin{equation*}
\Delta _{f}\frac{R}{f}=\frac{R}{f^{2}}(2f-\frac{n}{2})-2\frac{\left\vert
\operatorname{Rc}\right\vert ^{2}}{f}-4\frac{\operatorname{Rc}(\nabla f,\nabla f)}{f^{2}}+2%
\frac{R|\nabla f|^{2}}{f^{3}}\leq -\frac{|\operatorname{Rc}|^{2}}{f}+4\frac{(1+\sqrt{%
n})^{2}}{f}.\vspace{-0.04in}
\end{equation*}%
If $\zeta \!\left( s\right) \!=\!s$ for $s\in \left[ 0,1\right] $, $\zeta
\!\left( s\right) \!=\!1$ for $s\in \left[ 1,\bar{s}-1\right] $, $\zeta
\!\left( s\right) \!=\!\bar{s}-s$ for $s\in \left[ \bar{s}-1,\bar{s}\right] $%
, then\vspace{-0.04in}%
\begin{equation*}
\frac{c}{2}\int_{0}^{\bar{s}}\zeta ^{2}\left( \frac{|\operatorname{Rc}|^{2}}{f}-4%
\frac{(1+\sqrt{n})^{2}}{f}\right) ds+\frac{1}{2}\int_{0}^{\bar{s}}\zeta
^{2}\left\vert S\right\vert ^{2}ds\leq 2n-\int_{0}^{\bar{s}}2\zeta \zeta
^{\prime }\left\langle \nabla f,S\right\rangle ds.\vspace{-0.04in}
\end{equation*}%
Let $\gamma \left( 0\right) =x$, $\gamma \left( \bar{s}\right) =y$, and $%
\bar{s}=d\left( x,y\right) $. Then $1-c\leq C\leq 1+c$; the lower by $\frac{R%
}{f}\leq 1$ and the upper since for a minimal geodesic $\bar{\gamma}(s)$, $%
s\in \lbrack 0,\bar{s}]$, from $x$ and $y$,\vspace{-0.04in}%
\begin{equation*}
C\bar{s}\leq \int_{0}^{\bar{s}}\left( \left\vert \gamma ^{\prime }\left(
s\right) \right\vert ^{2}+c\frac{R(\gamma (s))}{f(\gamma (s))}\right) ds\leq
\int_{0}^{\bar{s}}\left( \left\vert \bar{\gamma}^{\prime }\left( s\right)
\right\vert ^{2}+c\frac{R(\bar{\gamma}(s))}{f(\bar{\gamma}(s))}\right)
ds\leq (1+c)\bar{s}.\vspace{-0.04in}
\end{equation*}%
Let $f\left( O\right) =\min_{\mathcal{M}}f\leq \frac{n}{2}$ and $r=d\left(
\cdot ,O\right) $. Then $\left\vert \nabla f\right\vert \left( z\right) \leq
\sqrt{f\left( z\right) }\leq \sqrt{\frac{n}{2}}+r\left( z\right) $. Since $%
\left\vert S\right\vert \leq \sqrt{C+c}$ and $r\left( \gamma \left( s\right)
\right) \leq \min \{r\left( x\right) +s\sqrt{C+c},r\left( y\right) +\left(
\bar{s}-s\right) \sqrt{C+c}\}$,\vspace{-0.04in}%
\begin{eqnarray*}
-\int_{0}^{\bar{s}}\!\zeta \zeta ^{\prime }\left\langle \nabla
f,S\right\rangle ds\!\! &\leq &\!\!\int_{0}^{1}\!s\sqrt{f\left( \gamma
\left( s\right) \right) }\left\vert S\left( s\right) \right\vert ds+\int_{%
\bar{s}-1}^{\bar{s}}\!\left( \bar{s}-s\right) \sqrt{f\left( \gamma \left(
s\right) \right) }\left\vert S\left( s\right) \right\vert ds \\
&\leq &\!\!\tfrac{1}{2}\sqrt{C+c}\left( \sqrt{2n}+r\left( x\right) +r\left(
y\right) +2\sqrt{C+c}\right) .
\end{eqnarray*}%
Let $A=\sqrt{C+c}$. Since $f\left( \gamma \left( s\right) \right) \geq
f\left( O\right) $ and $\bar{s}=d\left( x,y\right) $, we have\vspace{-0.04in}%
\begin{equation*}
\int_{0}^{\bar{s}}\frac{\zeta ^{2}|\operatorname{Rc}|^{2}}{f}ds\leq \frac{4(1+\sqrt{n%
})^{2}d\left( x,y\right) }{f\left( O\right) }+\frac{4(\sqrt{n}+A)^{2}}{c}+%
\frac{2A\left( r\left( x\right) +r\left( y\right) \right) }{c}.\vspace{%
-0.04in}
\end{equation*}%
Take $x=O$ and $\bar{s}=r\left( y\right) \geq 2\sqrt{\frac{n}{2}}$. Then $%
d\left( \gamma \left( s\right) ,y\right) \leq \frac{r\left( y\right) }{2}$
for $s\in \lbrack \frac{2A-1}{2A}\bar{s},\bar{s}]$ and\vspace{-0.04in}%
\begin{equation*}
\frac{(\frac{r\left( y\right) }{2A}-1)\min\limits_{s\in \lbrack (1-\frac{1}{%
2A})\bar{s},\bar{s}]}|\operatorname{Rc}|^{2}\left( \gamma \left( s\right) \right) }{(%
\sqrt{\frac{n}{2}}+\frac{3r\left( y\right) }{2})^{2}}\leq \int_{(1-\frac{1}{%
2A})\bar{s}}^{\bar{s}-1}\frac{|\operatorname{Rc}|^{2}\left( \gamma \left( s\right)
\right) }{f\left( \gamma \left( s\right) \right) }ds\leq \operatorname{Const}\left(
r\left( y\right) +1\right) .\vspace{-0.04in}
\end{equation*}%
Thus there exists $C<\infty $ such that for any $y\in \mathcal{M}$ with $%
r\left( y\right) \geq \max \{\sqrt{2n},3A\}$, there exists a point $z\in
\mathcal{M}$ with $d\left( z,y\right) \leq \frac{r\left( y\right) }{2}$ and $%
|\operatorname{Rc}|\left( z\right) \leq C\left( r\left( y\right) +1\right) $.%
\footnotetext[1]{%
Address. Bennett Chow, Bo Yang: Math. Dept., UCSD; Peng Lu: Math. Dept., U
of Oregon.}

\end{document}